\documentstyle[12pt]{article}
\begin{document}
     
\newtheorem{th}{Theorem}
\newtheorem{cor}{Corollary}
\newtheorem{prop}{Proposition}
\newtheorem{lem}{Lemma}
     
\newcommand{\p}{\partial}
\newcommand{\pmi}{\partial^{-1}}

\title{On a ring of formal pseudo-differential operators}
\thanks{This work was supported by the INTAS, EPSRC and RFFI grants.}
\author{A. N. Parshin}
\date{Steklov Mathematical Institute \\
Russian Academy of Sciences\\
     Gubkina str., 8\\
     Moscow, GSP-1 \\
        Russia}
\maketitle

We begin with a definition of the higher local fields arising in
algebraic geometry as purely local objects attached to algebraic
varieties of arbitrary dimension.
                  
\par\smallskip

 Let $K$ and $k$ be some fields. We say that $K$ has a structure
 of a {\it $n$-dimensional local field} with the (last) residue field $k$
if either $n=0$
 and $K=k$ or $n\geq 1$ and $K$ is the quotient field of a complete discrete
 valuation ring ${\cal O}_K$ whose residue field $\bar K$ is a local field
 of dimension $n-1$ with the last residue field $k$.

\par\smallskip
A typical example is the field of iterated Laurent power series
$$ K = k((x_1))\dots ((x_n)) $$
with the following local structure:
$$ {\cal O}_{K} = k((x_1))\dots ((x_{n-1}))[[x_n]] $$
$$\bar K = k((x_1))\dots ((x_{n-1})) $$
$$ \dots \dots \dots$$

If the characteristics of all the residue fields are equal, then the
     field $K$ must be of this form (for other examples and a more general
classification theorem see \cite[ch.2]{FP}).

This construction has first appeared in algebraic number theory for
     $n = 1$ and  it was later used in the theory of algebraic curves over
an arbitrary field. It can be defined for varieties and schemes of
     arbitrary dimension (see a survey \cite{FP}) and has  numerous
     applications to the problems of algebraic geometry, both arithmetical
     and geometrical.
     
     From this viewpoint, it would be reasonable to restrict ourselves
     by the commutative fields in the definition given above.
     Nevertheless, in the class field theory we meet the rings of this kind but
     the non-commutative ones. We mean the skew-fields which are
     finite-dimensional over their center $K$.
     Thus the $K$ will be a (commutative) local field and the skew-fields
     will represent the elements of the Brauer group of  the field $K$
     (see \cite{P}).
     
     The main purpose of this note is to point out onto another class
of non-commutative  local  fields  arising  in  the  theory of differential
equations and to show that these skew-fields hold many features of the
     commutative  fields.  We  define  a   skew-field   $P$   of   formal
pseudo-differential operators in $n$  variables which will be an
$2n$-dimensional  local  field. In contrast with the previous examples, it
will be infinite-dimensional over its center.
     The main properties of the skew-field $P$ and of it's "order"$E$
    are given in \S\S 1 and 2 \footnote{There are other approaches to
the rings of pseudo-differential operators in $n$
variables, see \cite{GQS}, \cite{GS} ¨ \cite{D}.
     We have used some arguments from \cite{GS}. The work
     \cite{D} became known to the author after the present text had been
written.}.
     
     We note that for the case of one variable a wider class of twisted
     pseudo-differential operators was considered in general algebra
     \cite{C} but there the purposes and motivations are purely algebraic
     and have no relations neither to algebraic geometry, nor to
differential equations.
     
     It  would  be  very  interesting  to   get   a   classification   of
non-commutative
     local fields at least for the case of equal characteristics of the
     residue fields \footnote{Several interesting results along this line
were  recently  proved by A. B. Zheglov, see his papers in Uspekhi
Matem. Nauk, {\bf 54}(1999), N 4, pp. 177-178 and Izvestija RAN (to appear)}. The commutative case is exhausted by the
     example  of  the  Laurent  power  series. This classification should
contain  as  a  particular  case  a   description   of   the   rings   of
pseudo-differential   operators  and  an  explicit  construction  of
skew-fields from the Brauer group.
     
     Another interesting problem is to find a generalization of the
construction of  commutative local field from a chain of subvarieties
     of an algebraic variety \cite{FP}.
     It would be reasonable to define a class of non-commutative schemes
     and to extract  the non-commutative local fields by an appropriate
localization process.
     
     The paper is completed by  \S \S 3 and 4,  where we show that
     the Kadomtsev-Petviashvili hierarchy (the dynamical system defined
     on the space $\p + E_{-}$ for the case $n = 1$) can be extended
to the space $P^n$. The usual properties of  the  KP  hierarchy  will  be
preserved under the new circumstances, in particular, the existence of
infinitely many conservation laws, the zero curvature presentation as the
Zakharov-Shabat equations, etc. A connection of this generalization of the
     KP hierarchy and of some natural Poisson structures on the space
$P^n$ is discussed in \S 4.
     
     A part of this work was done during my stay at the International
     Center of Theoretical Physics (Trieste, Italy) in September-October
     1997. I am very much grateful to M. S. Narasimhan for the kind
     invitation and hospitality.
     
     \section{Skew-field $P$ and ring $E$}
     
     Let $A$ be an associative not necessary commutative ring and let
     $d: A \rightarrow A$ be its derivation. Let us first introduce the ring
     $A((\pmi))$ of formal pseudo-differential operators with
     coefficients from $A$  as a left $A$-module of all formal
     expressions
     $$ L = \sum_{i > -\infty}^{n} a_i \p^i,~a_i \in A. $$
     Then a multiplication can be defined according to the Leibnitz rule:
     $$ (\sum_i a_i\p^i)(\sum_j b_j\p^j) =
     \sum_{i,j;k \ge 0} {i\choose k}a_i d^k(b_j)\p^{i + j - k}.$$
     Here we put
     $${i\choose k} = \frac{i(i-1)\dots (i-k+1)}{k(k-1)\dots 1},
     ~\mbox{if}~k > 0,~{i\choose 0} = 1.$$
     Particularly, for  $a \in A$:
     $$[\p,a] = \p a - a\p = d(a),$$
     (Heisenberg commutative relation),
     $$[\p^{-1},a] = \p^{-1} a - a \p^{-1} = -d(a)\p^{-2} + d^2(a)\p^{-3}
     - \dots .$$
     It can be checked that $A((\p^{-1}))$ will be again an associative
     ring.
     For  the  case  of  a  ring $A$ of functions, this ring was introduced
     by I. Schur \cite{Sch}. Later the ring $A((\p^{-1}))$
     was many times rediscovered and studied (see \cite{GD}, \cite{SS},
     \cite{Mum}, \cite{Mul}).
     The case of a non-commutative ring $A$
     was considered at the beginning of  the  paper  \cite{Mul}.
     
     We can iterate the constuction and this gives an opportunity to
     consider $n$ variables and to introduce the following rings.
     
     {\sc Definition 1}. Let $k$ be a field and  let $x_1,\dots x_n$ be
     some commuting variables. We put
     
     $$ P =  k((x_1))\dots ((x_n))((\p_1^{-1}))\dots ((\p_n^{-1})),$$
     $$ E =  k[[x_1, \dots x_n]]((\p_1^{-1}))\dots ((\p_n^{-1})). $$
     
     Then $E$ is a subring in $P$. Here $k((x_1))\dots ((x_n))$ is the
     ring of iterated Laurent power series and $k[[x_1, \dots x_n]]$ is
     the ring of Taylor power series. In the first case the order of
     variables $x_1, \dots x_n$ is significant.

     If $ L = \sum_{i \le m} a_i \p_n^i$ and
     $a_m \ne 0$, then $m = \mbox{ord}(L)$ will be called the {\it order}
     of the operator $L$. The function ord(.) defines  a decreasing
     {\it filtration} $P_{.}: \dots \subset P_{-1} \subset P_0 \subset \dots$
     of vector subspaces $P_i = \{L \in P~: \mbox{ord}(L) \le
i \} \subset P$.
     
     We have  the decomposition of $P$ and, accordingly of $E$, in a direct
sum of
     subspaces
$$P = P_{+} + P_{-}, $$
$$ E = E_{+} + E_{-},~E_{\pm} = E \cap P_{\pm},$$
where $P_{-} = \{L \in P: \mbox{ord}(L) < 0 \}$~and~$P_{+}$~consists of the
operators containing only $ \ge 0$~powers~of~$\p_n$.

Let us define the {\it highest term} of an operator $L$
by induction on $n$.
If $ L = \sum_{i \le m} a_i \p_n^i$~and~$\mbox{ord}(L) = m$ then  we
put
     $$\mbox{highest term}(L) = (\mbox{highest term}(a_m))\p_n^m.$$
     
     The highest term has the following form $f\p_1^{m_1}\dots \p_n^{m_n},
~f \in k((x_1))\dots ((x_n)),$ $~f \ne 0$ and one can define
     $$ \nu (L) = (m_1, \dots m_n) \in \mbox{{\bf Z}}^n. $$
     
     \begin{prop}. The rings $P$ and  $E$  have the following properties:
     
     i)  $P$  is  an  associative  skew-field;  an  operator $L \in E$ is
invertible in the ring $E$
     $\Leftrightarrow$ the coefficient $f$ in the highest term of $L$ is
invertible in the ring $k[[x_1, \dots x_n]]$

     ii) $$\mbox{ord}(LM) = \mbox{ord}(L) + \mbox{ord}(M),$$
        $$ \mbox{ord}(L + M) \le \mbox{max}(\mbox{ord}(L), \mbox{ord}(M))$$

     iii) $$\nu(LM) = \nu(L) + \nu(M),$$
        $$ \nu(L + M) \le \mbox{max}(\nu(L), \nu(M))$$
     for the lexicographical order in $\mbox{{\bf Z}}^n$.

     iv) Let $\nu(L)$ be divisible by $m \in \mbox{{\bf N}}$ such that

$(m,\mbox{char}(k) = 1$. If the coefficient $f$ in the highest
     term of the operator $L \in P~(\mbox{or}~\in E)$ is a $m$-th power
     in $k((x_1))\dots ((x_n))$(or correspondingly in $k[[x_1,\dots x_n]]$),
     then there exists a unique, up to  multiplication by
a  $m$-th root of unity,
     operator $M \in P~(\mbox{or}~\in E)$ such that
      $L = M^m$.
     \end{prop}
     {\sc Proof}. For the case $n = 1$, these claims are well known
 (see \cite{Sch}, \cite{Mum}, \cite{GS}). The general case
     can be done by a simple induction over $n$. We only note that
     the proof of associativity for the ring $P$ can be given exactly
     as for $n = 1$ (\cite{Mum}[ch. III, \S 11]) using  lemma 1.
     
     It is important to note that inequality
     $$\mbox{ord}([L,M]) \ge \mbox{ord}(L) + \mbox{ord}(M) - 1,$$
     known for $n = 1$,
     is in general wrong for $n > 1$ (see, for example, equality
     (5) ¢ \S 2).
     
     \par\smallskip
     
     Thus $P_0$ is a discrete valuation ring in
     $P$ with  residue field $$k((x_1))\dots ((x_n))((\p_1^{-1}))\dots
((\p_{n-1}^{-1}))$$
and we can introduce a structure
    of a  $2n$-dimensional local field on $P$
by an induction
(considering also the skew-fields with a non-necessarily equal number of
variables and derivations)

     Now let $V = k((x_1))\dots ((x_n))((z_1))\dots ((z_n))$. This is
     a field.
     {\it The symbol} of an operator  $L = \sum a_{i_1 \dots i_n} \p_1^{i_1}
     \dots \p_n^{i_n}$ is the series $\sigma(L) = \sum a_{i_1 \dots i_n}
     z_1^{i_1}\dots z_n^{i_n}$. For $i,j = 1,\dots n$ there are linear
     maps
     $$ d_i,~D_j: V \rightarrow V, \mbox{where}~d_i = \p/\p x_i;~D_j =  \p/\p
     z_j. $$
     
     Let us introduce several notations which are standard for the theory
    of differential operators
     in $n$ variables. We denote by $\alpha, \beta,
     \gamma$ the elements from $\mbox{{\bf Z}}^n$. Thus we will use
the following abbreviations for  operators
     $$\p^{\alpha} = \p_1^{\alpha_1}\dots \p_n^{\alpha_n},
     ~d^{\alpha} = d_1^{\alpha_1}\dots d_n^{\alpha_n},
     ~D^{\alpha} = D_1^{\alpha_1}\dots D_n^{\alpha_n},$$
     variables $z^{\alpha} = z_1^{\alpha_1}\dots z_n^{\alpha_n}$~
     and coefficients
     $${\alpha \choose \beta} = {\alpha_1 \choose \beta_1}\dots {\alpha_n
     \choose \beta_n},
      ~\alpha ! = \alpha_1 ! \dots \alpha_n !,$$
     where $\alpha = (\alpha_1, \dots \alpha_n).$
     We also write $\alpha \ge 0$ if all the $\alpha_i \ge 0$.

     \begin{lem}. Let $F,~G \in V$. Put
     $$ F * G = \sum_{\alpha \ge 0} \frac{1}{\alpha !} D^{\alpha}F
     d^{\alpha}G.$$
     For operators $L,~M \in P$ we then have
     $$ \sigma(LM) = \sigma(L)*\sigma(M). $$
     \end{lem}
     
     {\sc   Proof.}   It is sufficient  to check the lemma for  $L  =
a\p^{\alpha},~M = b\p^{\beta}$. Then $\sigma(L) = az^{\alpha},~\sigma(M) =
     bz^{\beta}$ and we get
     $$\sigma(L)*\sigma(M) = \sum_{\gamma \ge 0} \frac{1}{\gamma !}
D^{\gamma}\sigma(L)
     d^{\gamma}\sigma(M) = \sum_{\gamma \ge 0} \frac{1}{\gamma !} a
     D^{\gamma}(z^{\alpha}) \p^{\gamma}(b)z^{\beta} =$$
     $$
     \sum_{\gamma \ge 0}  {\alpha \choose \gamma}a\p^{\gamma}(b)
     z^{\alpha + \beta - \gamma } = \sigma(LM).$$
     
     Now we go to the constructions in the ring $P$ related to duality.
     
     {\sc Definition 2.} Let $L \in P$ ¨ $L = \sum a_{i_1 \dots i_n;
     j_1\dots j_n} x_1^{i_1}\dots x_n^{i_n}\p_1^{j_1}\dots \p_n^{j_n}$.
     Then the {\it residue} of the operator $L$ is defined as
     $$ \mbox{res}_P (L) = a_{-1\dots -1} \in k.$$
     
     If $L, M \in P$, then we put
     $$<L, M> = \mbox{res}_P (LM).$$

     \begin{prop}. The residue $\mbox{res}_P$ and the pairing $<.,>$ have
     the following properties:
     
     i)  $\mbox{res}_P$~is a  linear form on $P$ and
     $$\mbox{res}_P(L) = \mbox{res}_V(\sigma (L)) $$
     
     ii) for any   $L,~M \in P$
     $$ \mbox{res}_P ([L, M]) = 0.$$
     
     iii) the pairing~$<.,.>$~is a nondegenerate bilinear form on $P$.
     
     \end{prop}
     
{\sc Proof.} i) This claim is obvious. ii) Let $K
     = k((x_1))\dots ((x_n))$. Then $\mbox{res}_P(L) = \mbox{res}_{K/k}
(\mbox{res}_{P/K}(L))$.
Let us show that if $L \in [P,~P]$,
then
$$ \mbox{res}_{P/K}(L) \in d_1(K) + \dots + d_n(K).$$
     This gives what we need.
     By  lemma  1  and property i), it is enough to check that we have
     \begin{equation}
     \label{*}
F*G - G*F \in d_1(V) + \dots + d_n(V) + D_1(V) + \dots + D_n(V)
\end{equation} in the field $V$ of symbols.
     Remembering the definition of $*$,  we can restrict ourselves by
considerating  the term
 $D^{\alpha}(F)d^{\alpha}(G) - D^{\alpha}(G)d^{\alpha}(F)$  from
the commutator in the left hand side of (\ref{*}). Let us
     write $f \sim g$ if $f-g$ belongs to the right hand side of (\ref{*}).
     Since  $d_i$ and $D_j$ are derivations, we have
     $$D^{\alpha}(F)d^{\alpha}(G) - D^{\alpha}(G)d^{\alpha}(F) =
     D^{\alpha}(F)d^{\alpha}(G) - d^{\alpha}(F)D^{\alpha}(G)  \sim $$
     $$
     FD^{\alpha}d^{\alpha}(G) - Fd^{\alpha}D^{\alpha}(G)  = 0,$$
     because $D^{\alpha}d^{\alpha} = d^{\alpha}D^{\alpha}$. The proof
     given here is an adaptation to our situation of the arguments
from \cite{GS}[IV.50].

     iii) follows from
     $$ <x^{\alpha}\p^{\beta},~x^{\alpha'}\p^{\beta'}> = $$
     $$
     \left\{
     \begin{array}{ll}
     1 & \mbox{if}~\alpha + \alpha' = \beta + \beta' = (-1,\dots, -1)\\
     0 & \mbox{if for some}~i~\alpha_i + \alpha'_i + 1 < 0
     ~\mbox{or}~\beta_i + \beta'_i + 1 < 0.\\
     \end{array}
     \right.
     $$

     {\sc Remark 1.} The statement iii) shows that the skew-field $P$  is
     auto-dual. Thus it has the most important property of the
$n$-dimensional local fields.
     
     {\sc Remark 2.} Although we consider in this paper a purely formal
situation let us mention that the constructions and results given above
can be mostly extended to the case of ring ${\cal R} = R((\p_1^{-1}))
     \dots  ((\p_n^{-1}))$  where  $R$~is~  the  ring   of   $C^{\infty}$
functions on
     $\mbox{{\bf R}}^n$ rapidly decreasing at infinity. Proposition 1
     can be reformulated word by word. Concerning  proposition 2, one has
to consider the linear form
     $$tr: {\cal R} \rightarrow \mbox{{\bf R}}, $$
     $$ tr(\sum a_{i_1\dots i_n} \p_1^{i_1}\dots \p_n^{i_n}) =
     \int\limits_{\mbox{{\bf R}}^n} a_{-1,\dots -1} dx_1\dots dx_n$$
     instead of the operation $\mbox{res}_P$.
     
\section{The conjugacy theorems}
     We assume that the characteristic of the field $k$ is equal to 0.
Denote by $(<m)$ an arbitrary operator of an order $<m$.
     \begin{th}.
     Let $L_1 \in \p_1 + E_{-}, \dots  L_n \in \p_n + E_{-}$
     or $L_1 \in \p_1 + {\cal R}_{-}, \dots  L_n \in \p_n + {\cal R}_{-}$.
     Then the following conditions are equivalent:
     
     i) $\mbox{for any}~i,j~ [L_i, L_j] = 0$
     
     ii)  There exists an operator  $S \in 1 + E_{-}$
     (or $\in 1 + {\cal R}_{-}$) such that
     $$L_1 = S^{-1}\p_1 S, \dots L_n = S^{-1}\p_n S.$$
     
     If $S,S'$ are two operators from ii) then
     $$ S'S^{-1} \in 1 + k((\p_1^{-1}))\dots ((\p_n^{-1})) \cap E_{-}$$
     in the first case and
     $$ S'S^{-1} \in 1 + \mbox{{\bf R}}((\p_1^{-1}))\dots ((\p_n^{-1}))
     \cap {\cal R}_{-}$$ in the second case.
     \end{th}
     
       {\sc Proof} can be done by the same reasoning in both cases and,
to be concrete, we take the ring $E$. Denote by
     $R$ the ring $k[[x_1, \dots x_n]]$~and by  $E'$ the ring $R((\p_1^{-1}))
     \dots ((\p_{n-1}^{-1}))$.
     
     Now let
     \begin{equation}
     \label{1}
 L_i = \p_i + M_i = \p_i + a_i\p_n^m + (<m),~m < 0,~
     a_i \in E'.
     \end{equation}
A straightforward computation gives
     \begin{equation}
     \label{2}
[L_i, L_j] = (\frac{\p a_j}{\p x_i} -\frac{\p a_i}{\p x_j})\p_n^m +
    (<m),
     \end{equation}
     where we have extended the derivations  $\frac{\p}{\p x_i}$ from the ring
      $R$ to $E$ as above.
     
     If
\begin{equation}
     \label{3}
S = 1 - P,~P = b\p_n^m,~b \in E',
     \end{equation}

 then it is easy to check that
     \begin{equation}
     \label{4}
[\p_i, P] = \frac{\p b}{\p x_i}\p_n^m.
     \end{equation}
     Now we are ready to prove the theorem using  subsequent
approximations in powers of~$\p_n$.
     
     Let the operators $L_i$ satisfy the condition i) and have the form
as in (\ref{1}).
     We will look for an~$S$ of the form (\ref{3}). We have
     $$S^{-1}L_i S = (1 + P + P^2 + \dots)L_i(1 - P) = $$
     $$L_i - [L_i, P] - P[L_i, P] - P^2[L_i, P] - \dots = $$
     $$\p_i + M_i - [\p_i, P] - [M_i, P] - P[\p_i, P] - \dots .$$
     According  to  proposition  1,  all  the terms, except the first three
, belong to
     $E_{<m}$. Hence, we get
     $$S^{-1}L_i  S  =  \p_i  +  (a_i - \frac{\p b}{\p x_i})\p_n^m +
     (<m).$$
     Since the  commutators from (\ref{2}) are all equal to zero, it implies
     that there exists   $b \in E'$ such that
     $$ a_i = \frac{\p b}{\p x_i}, \mbox{for all}~i = 1, \dots n.$$
     It means that if we choose $S$ with such $b$, then the second term  in
our decomposition is equal to zero and thus the operators
     $L_i$ are conjugate to $\p_i$ up to terms of the next order.
     
The operators $S$ obtained in this way  step by step can be multiplied
inside the group
     $1 + E_{-}$ and the result is a solution of our problem. The inverse
implication $ii) \Rightarrow i)$ is obvious.
     
     In order to get the second statement of the theorem it is enough to
note that the conditions imply
     $$[S'S^{-1}, \p_i] = 0,\mbox{for all}~i = 1, \dots n.$$
     and, again applying  the formula (\ref{4}) by induction, we get that
     all coefficients of the operator
      $S'S^{-1}$ are  constant.
     
     This completes our  proof.

     Now denote by $Z(M)$ the ring of operators commuting with every
     operator from $M$.
     \begin{cor}.
Let $L_1 \in \p_1 + E_{-}, \dots  L_n \in \p_n + E_{-}$ and the condition i)
from theorem 1 is satisfied. Then the ring
     $Z(L_1, \dots L_n) \subset E$ is commutative and is equal to  $k((L_1^{-1}))
     \dots ((L_n^{-1}))$.
     \end{cor}
     {\sc Proof.} By theorem 1, we can assume that $L_i =
     \p_i$~for all~$i$. It is obvious that $Z(L_1, \dots L_n) \supset
     k((\p_1^{-1}))\dots ((\p_n^{-1}))$. The equality of these rings will
     follow from the next refinement of the formula (\ref{4}).

     Let
     $$ P = b\p_n^m + (<m),~b \in E',~m \in \mbox{{\bf Z}}.$$
     Then
     $$[\p_i, P] = \frac{\p b}{\p x_i}\p_n^m + (<m).$$
     
     \begin{prop}. The center of the ring $P$ is equal to $k$.
     \end{prop}
     
     {\sc Proof} is done using the formula just given
     (it is also valid in $P$). Namely, if   $L  \in
Z(P)$, then $L$ commutes with all the $\p_i$ and, therefore, it has  constant
coefficients. Because $L$ does also commute with all $x_i$, we get that
 $L \in k$.
     
    \par\smallskip
     
     For the case $n = 1$,  theorem 1 was proved by M. Sato \cite{SS},
\cite{Mul}.
     In this form, the result cannot be extended to the skew-field $P$.
Already for $n = 1$,  there are non-trivial obstacles  to the conjugacy
of  operators, which do not exist for the case of regular coefficients.
In  Sato's paper there is an extension of  theorem 1 for $n = 1$,
which uses so-called quasi-regular operators
(both for $L$ and $S$). Here we give another simple result of this kind
but also for $n = 1$.

     \begin{th}. Let $L, M$  be  operators of orders
     $\ne 0$ from the ring  $P$ of one variable.
     Assume that they have the form  $ \p^m + (< m-1)$.
Then the following conditions are equivalent:
     
     i) There is $S \in 1 + P_{-}$ such that $M = S^{-1}LS$
     
     ii) $\mbox{ord}(L) = \mbox{ord}(M) =: k$ and for all $m \in
     \frac{1}{k}\mbox{{\bf Z}}$
     $$ \mbox{res}_P(L^m) = \mbox{res}_P(M^m).$$
     \end{th}
     
     {\sc Proof.} The implication $i) \Rightarrow ii)$ follows from
     propositions 1 and 2 of section 1. To get the back implication, we
     reduce our claim to the case $k = 1$, taking the root
     (proposition 1, iv)),
.
     Next, we proceed as in the proof of theorem 1 (we only put $x_1 = x$
     ~and~$\p_1 = \p$).
     
     Let
     $$ L = M + a\p^m +(<m),~m \le -1.$$
     Then
     $$L^{-m} = M^{-m} + a\p^mM^{-m-1} + Ma\p^mM^{-m-2} + \dots
     + M^{-m-1}a\p^m + (<m) =$$
     $$ma\p^{-1} + M^{-m} +  (<m),$$
     because $M = \p + (<0)$. The conditions of the theorem show that
     $$\mbox{res}_{K}(a) = 0~\mbox{in the field}~K = k((x)).$$
     Therefore, there exists an element $b \in K$ such that $\p b/\p x = a$.
     
     Looking for $S$ as $1 + b\p^m$, we get as above
     $$S^{-1}LS = M + (a - \p b/\p x)\p^m + (<m) =
     M + (<m).$$

     It is interesting to compare this result with the classification of the
     conjugate    elements   in   division   algebras   $D$   which   are
finite-dimensional over their center $K$. We have
     
     \begin{th}. Let $X, Y \in D^*$. Then the  following  conditions  are
equivalent:
     
     i) There exists  $U \in D^*$ such that
     $$X = U^{-1}YU$$
     
     ii) $[K(X):K] = [K(Y):K] $ and for all $i \in \mbox{{\bf N}}$
     $$ \mbox{Tr}_{D/K}(X^i) = \mbox{Tr}_{D/K}(Y^i).$$
     \end{th}
      
     {\sc Proof.} The conditions in ii) mean that the elements  $X$ and
     $Y$ have the same minimal polynomials over $K$ and thus
     are conjugate in an algebraic closure of the field $K$. By the
Skolem-Noether  theorem \cite{B}[ch. VIII], there exists an element $U$ of
the algebra $D$ which realize this conjugation inside the $D$.

     {\sc Remark 3.} We immediately see that the conjugacy conditions in
 theorems  2  and  3  have  something  common.  The  valuation ord(.)
corresponds to the degree
     [:] and the residue $\mbox{res}_P$ corresponds to the trace
$\mbox{Tr}_{D/K}$.
     To understand an origin of this analogy, one has to remember
the adelic construction of the  direct image for symbols and differentials
for the morphisms of schemes of relative dimension $k$\cite{FP}.
     The degree [:] is a dimension and ord(.) is a valuation
and they are connected by a boundary map from K-theory.
This map from $\mbox{K}_1$ to $\mbox{K}_0$ will appear when one defines
the direct image of the K-sheaves for $k = 1$. In its turn, the trace and
residue is used for a definition of the direct image of differential forms,
     correspondingly, for $k = 0$~and~$k \ge 1$.
     (The case of the map of a surface onto a curve was considered  in
     \cite{O}).

At last, we note that the relative dimension (in the sense of
     local fields)
of the skew-field $D$ over its center $K$  is  equal  to  zero  and  this
dimension for the skew-field $P$ over its center $k$ is greater then zero.
     
\section{A generalization of the KP hierarchy}
                                 
     The vector space  $P$ and its power $P^n$ can be considered as
     infinite-dimensional varieties over the field $k$. Particularly, in
     every point $L = (L_1, \dots ,L_n) \in P^n$, the tangent space $T_{L}$
     is again $P^n$. The variety $P^n$
     contains a subvariety
     $$ P' = \{L \in P^n : \mbox{for all}~i,~j~[L_i, L_j] = 0 \}. $$
     Let us denote by $m$ (and also by $k$ and $l$) the multi-index
     $(m_1, \dots ,m_n)$ with nonnegative $m_i$ and
     introduce on $P^n$ the vector field $V^{m}$.
     
     {\sc Definition 3.}
     $$V^{m}_{L} = ([(L_1^{m_1}\cdots L_n^{m_n})_{+}, L_1],
\dots ,[(L_1^{m_1}\cdots L_n^{m_n})_{+}, L_n]) \in T_{L}. $$
     This field defines a dynamical system on $P^n$
which has the following form
     \begin{equation}
     \label{6}
 \frac{\p L}{\p t_{m}}  = V_{L}^{m}.
     \end{equation}
     
If we want to speak about the solutions
of the system (\ref{6}) in our purely formal setting,
we have to consider the operators ${L}$ as belonging to the extended
phase space $P^n \otimes k[[\dots, t_{m},
     \dots]]$. This space contains an infinitely many "times" $t_{m}$.
     We will also use an abbreviation
     $L^{m} =  L_1^{m_1}\cdots L_n^{m_n}$.
     
     \begin{prop}. Let $L \in P'$. Then we have

i)  if $L$ satisfies to the system (\ref{6}) then $L$ also satisfies
     to the system
     \begin{equation}
\label{7}
 \frac{\p (L^m)_{+}}{\p t_{k}} -
 \frac{\p (L^k)_{+}}{\p t_{m}} = [( L^k)_{+},
     (L^m)_{+}]
\end{equation}
     
ii) if $S \in 1 + P_{-}$ satisfies to the equation
     \begin{equation}
     \label{8}
 \frac{\p  S}{\p t_ m}    =    -(S\p_1^{m_1}\cdots
\p_n^{m_n}S^{-1})_{-}S,
\end{equation}
then $L = (S\p_1S^{-1}, \dots ,S\p_nS^{-1})$ satisfies to the system
      (\ref{6}).
     \end{prop}
     
{\sc Proof.} i). Since both parts in (\ref{6}) are derivations,
we have
 $$\frac{\p L^m}{\p t_{k}} = [(L^k)_{+},  L^m] $$
and simultaneously
$$\frac{\p (L^m)_{+}}{\p t_k} = [(L^k)_{+},
     L^m]_{+}. $$
     Thus
 $$\frac{\p (L^m)_{+}}{\p t_k} -
 \frac{\p (L^k)_{+}}{\p t_{m}}
- [(L^k)_{+}, (L^m)_{+}] = $$
$$([(L^k)_{+}, L^m] -
[(L^m)_{+}, L^k]
- [(L^k)_{+}, (L^m)_{+}])_{+} = $$
 $$ ([(L^m)_{+} -  L^m, (L^k)_{+}
-  L^k])_{+} =
 [(L^m)_{-}, (L^k)_{-}]_{+} = 0,$$
     because for all $L \in P'$ we have $[L^m, L^k] = 0$.
     
    ii). Since
     
 $$\frac{\p  S}{\p t_m} = -(L_1^{m_1}\cdots L_n^{m_n})_{-}S =
- (L^m)_{-}S,$$
     we get
 $$\frac{\p  L_{i}}{\p  t_{m}}  =  \frac{\p S}{\p t_{m}}\p_i
S^{-1} -  S\p_i S^{-1}\frac{\p S}{\p t_{m}} S^{-1} = $$
     $$ - [(L^m)_{-}, L_i] =  [(L_1^{m_1}\cdots L_n^{m_n})_{+}, L_i].$$
     
     \begin{prop}. Let $L \in P'$. Then
     
     i)  $V^{m}_{L} \in T_{P',L},$
     
     ii) vector fields $V^{m}\mid_{P'}$  and $V^n\mid_{P'}$ commute.
     
\end{prop}
     
     {\sc Proof.} i). If  $L \in P'$, then the tangency condition
     for the vector $(X_1, \dots X_n)$ at the point $L$ means that for
     all pairs $i,j$ we have
     $$   [L_i + \epsilon X_i, L_j + \epsilon X_j] = 0,~\mbox{up to}
     ~\epsilon^2,$$
     and thus $[X_i, L_j] = [X_j, L_i]$. The Jacobi identity for the Lie
algebra $P$ shows that any field of the form $([U, L_1], \dots ,[U, L_n])$
will be tangent to $P'$.
     
     ii). Let $L(\dots ,t_{m},\dots)$ be a formal solution of
      (\ref{6}) such that $L(0) \in P'$. Then
     $$ \frac{\p}{\p t_m}(\frac{\p}{\p t_k}L_i) -
      \frac{\p}{\p t_k}(\frac{\p}{\p t_m}L_i) =
      \frac{\p}{\p t_m}[(L^k)_{+}, L_i] -
        \frac{\p}{\p t_k}[(L^{m})_{+}, L_i] =$$
     $$[\frac{\p}{\p t_m}(L^k)_{+}, L_i]  +
       [(L^k)_{+}, \frac{\p}{\p t_m}L_i]
     -  [\frac{\p}{\p t_k}(L^{m})_{+}, L_i]
     - [(L^{m})_{+}, \frac{\p}{\p t_k}L_i]  = $$
     $$[\frac{\p}{\p t_m}(L^k)_{+} -
         \frac{\p}{\p t_k}(L^m)_{+}, L_i] +
       [(L^k)_{+}, [(L^{m})_{+}, L_i]] -
       [(L^{m})_{+}, [(L^k)_{+}, L_i]] =$$
     $$  [[(L^{m})_{+}, (L^k)_{+}], L_i]] +
       [(L^k)_{+}, [(L^m)_{+}, L_i]] -
       [(L^{m})_{+}, [(L^k)_{+}, L_i]] = 0.$$
     Here  we  have  used the Zakharov-Shabat equations (\ref{7}) and the
Jacobi identity in the Lie algebra $P$.
     The proof is completed.
     
     {\sc Definition 4.}~Let $L \in P^n$. Then we put
     $$ H_{k}(L) = \mbox{res}_P(L^k).$$
     
     \begin{prop}. The functions $H_{k}$ are the integrals of motion
     for the system (\ref{6}) on the manifold $P^n$.
     \end{prop}
     
     {\sc Proof.}
     $$\frac{\p}{\p t_m}H_{k} =
     \mbox{res}_P(\frac{\p}{\p t_m}L^k ) = $$
     $$ \mbox{res}_P((\frac{\p}{\p t_m}L_1^{k_1})L_2^{k_2}\dots
     L_n^{k_n} + \dots + L_1^{k_1}\dots L_{n-1}^{k_{n-1}}\frac{\p}{\p t_m}
     L_n^{k_n})  =$$
    $$ \mbox{res}_P([L^m, L_1^{k_1}]L_2^{k_2}\dots
     L_n^{k_n} + \dots + L_1^{k_1}\dots L_{n-1}^{k_{n-1}}
     [L^m, L_n^{k_n}])  =$$
     $$\mbox{res}_P(L^mL_1^{k_1})\dots L_n^{k_n} -
     L_1^{k_1}\dots L_{n}^{k_{n}}L^m)  = 0.$$
     
     The  space  $P^n$  has  many  submanifolds,  which   are   invariant
with respect to  the  system  (\ref{6}).  One  can  construct them using the
following result.
     
     \begin{lem}. Let $V(L) = L_{i_1}\cdots L_{i_k}$  be a monome in
operators $L_1,
     \dots , L_n$~and~$U \in P$. Then the vector field $F$ such that~$F_L =
     ([U_{+},L_1],\dots , [U_{+},L_n])$ is tangent to the variety
     $M_V = \{ (L_1,\dots ,L_n) \in P^n : V(L)_{-} = 0 \}$.
     \end{lem}
     {\sc Proof}. Let $F_i = [U_{+},L_i],~i = 1,\dots , n$.
     The tangency condition has the following form
     $$((L_{i_1} + \epsilon F_{i_1})\cdots (L_{i_k}
     + \epsilon F_{i_k}))_{-} = 0,
     ~\mbox{up to }~\epsilon^2.$$
     The left hand side is equal to
     $$V(L)  +  \epsilon
     \sum_{j=1}^k L_{i_1} \cdots  L_{i_{j-1}}F_{i_{j}}L_{i_{j+1}}  \cdots
     L_{i_k}
     + \epsilon^2 \dots ,$$
     and, by application of the operation $(.)_{-}$,  we get
     $$(\sum_{j=1}^k  L_{i_1} \cdots L_{i_{j-1}}U_{+}L_{i_{j}}L_{i_{j+1}}
     \cdots L_{i_k}  -
     \sum_{j=1}^k  L_{i_1} \cdots L_{i_{j-1}}L_{i_{j}}U_{+}L_{i_{j+1}}
     \cdots L_{i_k})_{-} = $$
     (here all the terms, except the first two, will be cancelled)
     $$(U_{+}L_{i_1}\cdots L_{i_k} - L_{i_1}\cdots L_{i_k}U_{+})_{-} =
     ([U_{+},V])_{-} = [U_{+},V_{+}]_{-} = 0.$$
     
      This lemma gives an opportunity to define the different versions
of the KdV hierarchy for $n > 1$.
      
     {\sc Remark 4.} All definitions and results of this section
     will be valid for the system (\ref{6}) on the space ${\cal R}^n$
     (instead of $P^n$) and with~$\mbox{tr}(L^k)$~(instead of
      $\mbox{res}_P(L^k)$).

     If $n = 1$, then the system (\ref{6}) is the Kadomtsev-Petviashvili
hierarchy extended to the whole space $P$
     (it is common to consider the system on the affine space $\p + E_{-}$).
     Accordingly, (\ref{7}) and  (\ref{8})  are  the  Zakharov-Shabat  and
Sato-Wilson equations.
The last ones give rise to the system (\ref{6}) if we restrict it
to $((\p_1,\dots ,\p_n) +  E_{-}^n)
     \cap P'$ and apply theorem 1.

     \section{Poisson structures on $P^n$ and Hamiltonians $H_k$}

     We want to define a Poisson structure on the space $P^n$ and
let us first introduce
      the space ${\cal F}(P^n)$ of functionals on $P^n$
     as  the  vector  space  of  polynomial  functions  in (any  number)
coefficients
     (close to~$x_1^{i_1}\cdots x_n^{i_n}\p_1^{j_1}\cdots \p_n^{j_n}$)
     of operators $L_1, \dots , L_n \in P$. The bilinear form $<.,.>$
     on $P$ can be naturally extended to the $P^n$:
     $$<L,~M> = \sum_{i=1}^n <L_i,M_i>, $$
     where $L = (L_1,\dots, L_n),~M = (M_1,\dots,M_n) \in P^n$.
     
     If $F \in {\cal F}(P^n)$, then we can define the {\it gradient}
     $\nabla_F(L) \in T_L = P^n$ of the functional $F$ at a point $L \in P^n$
by the condition
     $$
     <M, \nabla_F(L)> = \frac{d}{d\epsilon} F(L + \epsilon M)
     \mid_{\epsilon = 0}, ~\mbox{for all}~M \in P^n.
     $$
     The structure of a Lie algebra on $P^n$ gives us an opportunity to
     define a {\it Poisson bracket}  $\{F,G\}  \in  {\cal  F}(P^n)$  for
all $F,~G \in
     {\cal F}(P^n)$:
     $$ \{F,G\}(L) = <L, [\nabla_F(L), \nabla_G(L)]>.$$

    The decomposition $P^n = (P^n)_{+} + (P^n)_{-}$ gives rise to a new
structure of a Lie algebra $[~,~]_R$ defined by an $R$-matrix (=
difference of the projectors
     $(P^n)_{+}$~and~$(P^n)_{-}$, see \cite{ST}). We have $[X,Y]_R =
     \frac{1}{2}[RX,Y] + \frac{1}{2}[X,RY] =
[X_{+},Y_{+}] - [X_{-},Y_{-}]$. Thus we have introduced a new Poisson structure
 $\{.,.\}_R$ on $P^n$.
     
     We need two simple facts.
     
     Let $A, B, C \in P^n$. Then
     \begin{equation}
     \label{9}
     <A, [B,C]> = <[A,B], C>.
     \end{equation}
     
     Indeed, in any ring $A[B,C] = [A,B]C + [B,AC]$ and
     we have to apply proposition 2, ii).

     \begin{lem}. If  $F_M \in  {\cal  F}(P^n)$ such that
     $F_M(L) = <L, M>$, then $\nabla_{F_M} = M$.
     \end{lem}
     {\sc Proof} immediately follows from the definition of  gradient
and non-degenerateness of the form $<.,.>$ on $P^n$.
     
     The functions $H_k$ from \S 3  belong to  ${\cal F}(P^n)$. Let us
remember that
     $k = (k_1, \dots ,k_n)$~and all $k_i \ge 0$.
     
     \begin{lem}. Let $L \in P'$. Then
     $$\nabla_{H_k}(L) = (U_1, \dots , U_n),$$
     $$\mbox{where}~
     U_i = k_iL_1^{k_1}\cdots
     L_{i-1}^{k_{i-1}}L_{i}^{k_i-1}L_{i+1}^{k_{i+1}}\cdots L_n^{k_n},
     i = 1,\dots, n. $$
     \end{lem}
     
     {\sc Proof.}
     $$<M,\nabla_{H_k}(L)> = $$
     $$ \frac{d}{d \epsilon}\mbox{res}_P
     ((L_1 + \epsilon M_1)^{k_1}\cdots (L_n + \epsilon M_n)^{k_n})
     \mid_{\epsilon = 0} =$$
     $$\mbox{res}_P(\sum_{i=1}^{n} L_1^{k_1}\cdots L_{i-1}^{k_{i-1}}
     (\sum_{j=0}^{k_{i}-1}
     L_{i}^{j}M_{i}L_{i}^{k_i-j-1})L_{i+1}^{k_{i+1}}\cdots L_n^{k_n}) =$$
     $$ \mbox{res}_P(\sum_{i=1}^n M_iL_1^{k_1}\cdots L_{i-1}^{k_{i-1}}
     \sum_{j=0}^{k_i-1}(L_i^{k_{i}-1})L_{i+1}^{k_{i+1}} \cdots
     L_{n}^{k_n}) =$$
     (using proposition 2, ii) and  that the  $L_i$  commute   with
each other)
     $$ \mbox{res}_P(\sum_{i=1}^n M_ik_iL_1^{k_1}\cdots L_{i-1}^{k_{i-1}}
     L_i^{k_{i}-1}L_{i+1}^{k_{i+1}} \cdots
     L_{n}^{k_n}) = $$
     $$ <M, U_i>.   $$
     
     Now we can compute the Poisson brackets for the functionals $H_k$.
     
     \begin{prop}. Let $L \in P'$. Then for all $k,l$~¨~$F \in
     {\cal F}(P^n)$
     $$i)~\{H_k,F\}(L) = 0,$$
     $$ii)~\{H_k, H_L\}_R(L) = 0.$$
     \end{prop}
     
     {\sc Proof.}
     $$i). \{H_k,F\}(L) = <L,[\nabla_{H_k}(L),\nabla_{F}(L)]> = $$
     (using (\ref{9}))
     $$ <[L,\nabla_{H_k}(L)],\nabla_F(L)> = 0,$$
     since by lemma 4 and by condition of the proposition
$[L,\nabla_{H_k}(L)] = 0$.
     $$ii).  \{H_k,H_l\}_R(L)  =   $$
     $$ \frac{1}{2}<L,[R\nabla_{H_k}(L),\nabla_{H_l}(L)]> +
     \frac{1}{2}<L,[\nabla_{H_k}(L),R\nabla_{H_l}(L)]> =$$
     $$ -\frac{1}{2}<[L,\nabla_{H_l}(L)],R\nabla_{H_k}(L)> +
     \frac{1}{2}<[L,\nabla_{H_k}(L)],R\nabla_{H_l}(L)> = 0,$$
     by the same argument as above.
     
     \begin{lem}. Let $H \in {\cal F}(P^n)$,~$\nabla_H(L) = U = (U_1,
     \dots , U_n)$~and, for all~$i = 1, \dots , n$, we have $[U_i, L_i] = 0$.
     Then the Hamiltonian system  $\frac{\p F}{\p t}  =  \{F,H\}_R,~F  \in
     {\cal F}(P^n)$  is of the form
     $$\frac{\p L_i}{\p t} = [(U_i)_{+}, L_i]$$
on the space  $P^n$.
     \end{lem}
     
     {\sc Proof.} Fix $M = (M_1, \dots , M_n) \in P^n$ and
     take the functional $F_M$ as $F$. Then
     $$\{F_M,H\}_R(L) = $$
     $$ <L,[\nabla_{F_M}(L),\nabla_{H}(L)]_R> =
     <L,[M,U]_R> = $$
     $$<L,[M_{+},U_{+}] - [M_{-},U_{-}]> = $$
     $$ - <[L,U_{+}],M_{+}> +
     <[L,U_{-}],M_{-}> =$$
     (using lemma 3 and (\ref{9}))
     $$- <[L,U_{+}],M_{+}> - <[L,U_{+}],M_{-}>  =$$
     (because by condition $[U,L] = 0$~¨~$U = U_{+} + U_{-}$)
     $$ = -<[L,U_{+}],M> = <[U_{+},L],M>.$$
     Since $\frac{\p F_M}{\p t} = <\frac{\p L}{\p t}, M>$ we get what we
     wanted.
     
\par\smallskip
     
     Combining together the lemmas 4 and 5, we get
     
     \begin{prop}. The Hamiltonian system
     $$\frac{\p F}{\p t_{m}} = \{F, H_m\}_R,~F \in {\cal F}(P^n)$$
     has the following form on the manifold $P'$ :
     $$\frac{\p L_i}{\p t_{m}} = m_i[(L_1^{}\cdots L_{i-1}^{m_{i-1}}
     L_i^{m_{i}-1}L_{i+1}^{m_{i+1}}
     \cdots L_{n}^{m_n})_{+}, L_i],~i= 1, \dots ,n.$$
     \end{prop}

     {\sc Remark 5.} This system coincides with the system (\ref{6}) \S 3
     only if $n = 1$. To present it in a Hamiltonian form, we have to
take  $n = 2$ and to consider only a part of the Hamiltonians.

     Namely, let now  $m,k,l \in \mbox{{\bf N}}$. We put
     $$H_m = \frac{1}{m}\sum_{k+l=m}{m \choose k}H_{k,l}.$$
     Then from lemma  4 we instantly find
     $$\nabla_{H_m}(L_1,L_2) =
     (\sum\limits_{i+j=m-1;i,j \ge 0} {m-1 \choose i}L_1^iL_2^j,
     \sum\limits_{i+j=m-1;i,j \ge 0}  {m-1 \choose j}L_1^iL_2^j)$$
     and this expression is of  the form $(U,U)$. Taking all together
we get the following result.
     \begin{prop}. The Hamiltonian system
     $$\frac{\p F}{\p t_{m}} = \{F, H_m\}_R,~F \in {\cal F}(P^2)$$
     has the following properties:
     
     i) The variety $P'$ is invariant and on $P'$  the system has the form:
     $$\frac{\p L_1}{\p t_m} = [(L_1 + L_2)^{m-1}_{+},L_1],$$
     $$  \frac{\p L_2}{\p t_m} = [(L_1 + L_2)^{m-1}_{+},L_2].$$
     ii) The functionals $H_m$ are the integrals of motion and
     $\{H_m,H_{m'}\}_R = 0$~on~$P'$.
     \end{prop}
     
     {\sc Remark 6.} Let us note that $P'$ is not a Poisson subvariety
     in $P^2$. This does not permit to present our dynamical system
     as a Hamiltonian one on the space $P'$. However, it is an open question
if there exists a Poisson structure
     on $P'$ such that the property ii) will be preserved and the system
would be Hamiltonian.

     {\sc Remark 7.} The decomposition
     $P^n = P_{+}^n + P_{-}^n$ used above is not a unique one.
     We note that the rings $P$~and~$E$ admit some non-standard
     decompositions (as vector spaces) into a direct sum
     of its subrings. For instance, for $n = 1$ we have,
besides $P = P_{+} + P_{-}$,
the following decompositions:
     $$P = xk[[x]]((\p^{-1})) + k[x^{-1}]((\p^{-1}))  =   $$
     $$ k[[x]]((\p^{-1})) + x^{-1}k[x^{-1}]((\p^{-1})),$$
   %  $$P = \{x^i\p^j : i \ge j \} +  \{x^i\p^j : i < j \} =
   %  \{x^i\p^j : i > j \} +  \{x^i\p^j : i \le j \}.$$
     $$P = xk[[x]](((x\p)^{-1})) + k[x^{-1}](((x\p)^{-1}))  = $$
     $$ k[[x]](((x\p)^{-1})) + x^{-1}k[x^{-1}](((x\p)^{-1})),$$
     
     It seems that the corresponding Hamiltonian systems have not yet
     been studied.


\begin{thebibliography}{30}

\bibitem{B}{Bourbaki  N ., {\it Alg\`ebre, Chapitre 8, Modules et anneaux
semi-simples}, Paris, Hermann, 1958}

\bibitem{C}{Cohn P. M., {\it Skew-fields}, Cambridge University Press, 1997}

\bibitem{D}{Dzhumadil'daev A. S., {\it Derivations and central extensions of
the Lie algebra of formal pseudo-differential operators}, Algebra and
     Analysis, 1994, v. 6, N 1, p. 140-158}
     
\bibitem{FP}{Fimmel T., Parshin A. N., {\it Introduction to the Higher Adelic
Theory}~(book in preparation)}

\bibitem{GD}{Gelfand I. M., Dickey L. A., {\it Fractional powers of operators
and Hamiltonian systems}, Function. Anal. Applic., 1976, v. 10, pp. 13-39}

\bibitem{GQS}{Guillemin V., Quillen D., Sternberg S., {\it The integrability
     of characteristics}, Commun. Pure Appl. Math., 1970, v. 23, p. 39-77}

\bibitem{GS}
{Guillemin V., Sternberg S., {\it Symplectic techniques in Physics}, Cambridge UP,
1984}

\bibitem{Mul}{Mulase M., {\it Solvability of the super KP equation and
a generalization of the Birkhoff decomposition}, Invent. Math.,
1988, vol. 92, pp. 1-46}

\bibitem{Mum}{Mumford D., {\it Tata lectures on Theta II}, Birkh\"auser,
Boston, 1984}
     
\bibitem{O}{Osipov D. V., {\it Adelic constructions of direct images of
differentials and symbols}, Matem. Sbornik, 1997, v. 188:5, pp. 59-84;
alg-geom/9802112}

\bibitem{P}{Parshin A. N., {\it Galois cohomologies and Brauer group of
local fields}, Proc. Steklov Math. Institute, 1990, v. 183, pp. 159-168}

\bibitem{ST}{Semenov-Tan-Shansky M. A., {\it What is the classical
$r$-matrix ?}, Function. Anal. Applic., 1983, v. 17, pp. 17-33}

\bibitem{SS}{Sato M., Sato Y.,~{\it Soliton Equations as Dynamical Systems
     on  Infinite  Dimensional  Grassmann Manifold}, Lecture Notes in Num.
Appl. Anal., 1982, vol. 5, pp. 259-271}

\bibitem{Sch}{Schur  I., {\it \"Uber  vertauschbare  lineare
Differentialausdr\"ucke}, Sitzungsber. der Berliner Math. Gesell., 1905,
     Bd. 4, S. 2-8}
\end{thebibliography}
\end{document}